\documentclass[reqno,A4,12 pt]{amsart}
\usepackage{amssymb}

\newcommand{\gs}[1]{\textbf{#1}}
\newcommand{\ds}{\displaystyle}
\newcommand{\fr}[2]{\frac{#1}{#2}}
\newcommand{\dfr}[2]{\dfrac{#1}{#2}}
\newcommand{\tfr}[2]{\tfrac{#1}{#2}}
\newcommand{\cd}{\cdot}
\newcommand{\cds}{\cdots}
\newcommand{\dsum}{\displaystyle \sum}
\newcommand{\dlim}[1]{\ds \lim_{#1}}
\newcommand{\ol}[1]{\overline{#1}}
\newcommand{\ul}[1]{\underline{#1}}
\newcommand{\simto}{\stackrel{\sim}{\to}}

\renewcommand{\l}{\left}
\renewcommand{\r}{\right}

\newcommand{\vsv}{\vspace{5mm}}
\newcommand{\vsb}{\vspace{2mm}}
\newcommand{\vsi}{\vspace{1mm}}
\newcommand{\q}{\quad}
\newcommand{\qq}{\qquad}
\newcommand{\sfr}[2]{\leavevmode\kern-.1em
  \raise.5ex\hbox{\the\scriptfont0 #1}\kern-.1em
  /\kern-.15em\lower.25ex\hbox{\the\scriptfont0 #2}}

\newcommand{\shf}{\sfr{1}{2}}
\DeclareMathOperator*{\dpi}{\prod}
\DeclareMathOperator*{\PI}{\textstyle \prod}

\newcommand{\pf}{\noindent {\bf Proof:}\q }
\newcommand{\supp}{\mathrm{supp}}
\newcommand{\del}{\partial}
\newcommand{\stab}{\mathrm{Stab}}
\newcommand{\pstab}{\mathrm{Stab}^{\mathrm{pt}}}
\newcommand{\h}{\mathfrak{h}}

\newcommand{\la}{\langle}
\newcommand{\ra}{\rangle}
\newcommand{\imp}{\Longrightarrow}
\newcommand{\rimp}{\Longleftarrow}
\newcommand{\abs}[1]{\lvert{#1}\rvert}
\newcommand{\any}{\forall}
\newcommand{\onetoone}{\longleftrightarrow}

\renewcommand\qedsymbol{\rule{3pt}{6pt}}
\renewcommand{\root}[1]{\sqrt{\mathstrut #1}}
\newcommand{\nthroot}[2]{{}^{#1}\! \! \! \root{#2}}

\DeclareMathOperator*{\ooplus}{\pmb{\oplus}}
\DeclareMathOperator*{\tensor}{\otimes}
\DeclareMathOperator*{\fusion}{\boxtimes}

\newcommand{\Z}{\mathbb{Z}}
\newcommand{\C}{\mathbb{C}}
\newcommand{\R}{\mathbb{R}}
\newcommand{\N}{\mathbb{N}}
\newcommand{\Q}{\mathbb{Q}}
\newcommand{\K}{\mathbb{K}}
\newcommand{\M}{\mathbb{M}}
\newcommand{\RM}{\mathrm{RM}}

\newcommand{\res}{\mathrm{Res}}
\newcommand{\End}{\mathrm{End}}
\newcommand{\vir}{\mathrm{Vir}}
\newcommand{\aut}{\mathrm{Aut}}
\newcommand{\wt}{\mathrm{wt}}
\newcommand{\tr}{\mathrm{tr}}
\renewcommand{\hom}{\mathrm{Hom}}
\newcommand{\ch}{\mathrm{ch}}
\newcommand{\SL}{\mathrm{SL}}
\newcommand{\id}{\mathrm{id}}
\newcommand{\ind}{\mathrm{Ind}}
\newcommand{\ad}{\mathrm{ad}}
\newcommand{\om}{\omega}
\newcommand{\be}{\beta}
\newcommand{\al}{\alpha}

\newcommand{\pii}{\pi \sqrt{-1}\, }
\newcommand{\hf}{\fr{1}{2}}
\newcommand{\Span}{\mathrm{Span}}

\newcommand{\w}{\omega}
\newcommand{\vacuum}{\mathrm{1\hspace{-3.2pt}l}}
\newcommand{\vac}{\vacuum}
\newcommand{\ising}{L(\tfr{1}{2},0)}
\newcommand{\irr}{\mathrm{Irr}}

\makeatletter \@addtoreset{equation}{section}
\renewcommand\theequation{\thesection.\arabic{equation}}

\theoremstyle{plain}
\newtheorem{thm}{Theorem}[section]
\newtheorem{prop}[thm]{Proposition}

\theoremstyle{definition}
\newtheorem{df}[thm]{Definition}

\theoremstyle{remark}
\newtheorem{rem}[thm]{Remark}

\renewcommand{\baselinestretch}{1.5}

\begin{document}

\title[]{A characterization of the moonshine vertex operator algebra by means of Virasoro frames
}

\author[C. H. Lam]{Ching Hung Lam $^\dagger$}%
\address[C.H. Lam]{Department of Mathematics, National Cheng Kung University,
Tainan, Taiwan 701}%
\email{chlam@mail.ncku.edu.tw}%

\author[H. Yamauchi]{  Hiroshi Yamauchi $^\ddag$}
\address[H. Yamauchi]{Graduate School of Mathematical
Sciences, The University of Tokyo, Tokyo, 153-8914, Japan}

\email{yamauchi@ms.u-tokyo.ac.jp}

\thanks{$^\dagger$ Partially
supported by NSC grant 95-2115-M-006-013-MY2 of Taiwan, R.O.C}%

\thanks{$^\ddag$ Supported by JSPS Research Fellowships for Young Scientists.}

\begin{abstract}
In this article, we show that a framed vertex operator algebra $V$
satisfying the conditions: (1) $V$ is holomorphic (i.e, $V$ is the
only irreducible $V$-module); (2) $V$ is of rank $24$; and (3)
$V_1=0$; is isomorphic to the moonshine vertex operator algebra
$V^\natural$ constructed by Frenkel-Lepowsky-Meurman \cite{flm}.
\end{abstract}
\maketitle

\section{Introduction}

The moonshine vertex operator algebra $V^{\natural}$ constructed by
Frenkel-Lepowsky-Meurman \cite{flm} is one of the most important
examples of vertex operator algebra (VOA). In the introduction of
their book\,\cite{flm}, Frenkel, Lepowsky and Meurman conjectured
that $V^{\natural}$ can be characterized by the following three
conditions:\\
(1) the VOA $V^{\natural}$ is the only irreducible ordinary module
for itself;\\
(2) the rank (or central charge) of $V^{\natural}$ is 24;\\
(3) $V^{\natural}_1=0.$

\medskip

While conditions (2) and (3) are clear from the construction,
condition (1) was proved by Dong \cite{d1} using the $48$ mutually
commuting Virasoro elements of central charge $1/2$ discovered in
\cite{dmz}. The discovery of the 48 commuting Virasoro algebras
inside $V^{\natural}$ also inspired the study of framed vertex
operator algebras \cite{dgh,m2}. In fact, a lot of progress have
been made on the study of the moonshine VOA using the 48 commuting
Virasoro algebras \cite{dgh,huang,lam3,ly,m3}.

In this article, we shall prove a weak version of FLM conjecture by
adding an extra assumption that $V$ contains $48$ mutually commuting
Virasoro elements of central charge $1/2$ (i.e., $V$ is a framed VOA).

\begin{thm}
Let $V$ be a framed vertex operator algebra satisfying the
conditions:\\ (1) $V$ is the only irreducible module for itself;
(2) $V$ is of rank $24$; and (3) $V_1=0$. Then $V$ is isomorphic
to the moonshine vertex operator algebra $V^\natural$.
\end{thm}

FLM's construction \cite{flm} of the moonshine VOA $V^\natural$ is
the first mathematical example of the so-called $\Z_2$-orbifold
construction and it is closely related to the lattice VOA associated
with the Leech lattice. Let $\Lambda$ be the Leech lattice and
$V_\Lambda$ the lattice associated with $\Lambda$. As shown in
\cite{flm}, we can extend the $-1$ map on $\Lambda$ to an
automorphism $\theta$ of $V_\Lambda$ by
\begin{equation}\label{eq:1.1}
  \theta: \al_1(i_1)\cdots \al_k{(i_k)}e^{\al} \mapsto
(-1)^{k}\al_1(i_1)\cdots \al_k(i_k)e^{-\al}
  \mbox{ for }\al_1,...,\al_k,\al\in \Lambda.
\end{equation}
Let $V_\Lambda^T$ be the unique irreducible $\theta$-twisted module.
Then the moonshine VOA $V^\natural$ is given by
\[
V^\natural= V_\Lambda^+\oplus V_\Lambda^{T,+},
\]
where $V_\Lambda^+$ and $V_\Lambda^{T,+}$ are the fixed point
subspaces of $\theta$ on $V_\Lambda$ and $V_\Lambda^{T}$,
respectively. By the construction, there is a natural involution
$t\in \aut(V^\natural)$ such that $t|_{V_\Lambda^+}=1$ and
$t|_{V_\Lambda^{T,+}}=-1$. This automorphism $t$ belongs to the $2B$
conjugacy class of the Monster group and the top weight of the
unique irreducible $t$-twisted module over $V^\natural$ is $1$
\cite{huang,lam3}. By performing the $t$-orbifold construction on
$V^\natural$, one can recover the Leech lattice VOA $V_\Lambda$
\cite{huang,lam3,ly}.

Our main strategy is to reverse the above construction and try to
obtain  the Leech lattice VOA $V_\Lambda$ by using the
$\Z_2$-orbifold construction on $V$. The key point is that for any
framed VOA $V$, one can easily define some involutions on $V$ by
using the frame structure. Such kind of involutions are often called
$\tau$-involutions or Miyamoto involutions (cf.\ \cite{m1}).
If $V_1=0$, we can define an
involution $\tau$ on $V$ such that the top weight of the irreducible
$\tau$-twisted module $V^T$ is $1$. Then by performing
$\tau$-orbifold construction on $V$, we obtain a VOA
\[
  V(\tau)=V^{\la \tau\ra}\oplus (V^T)^{\la \tau\ra},
\]
where $V^{\la \tau\ra}$ and $(V^T)^{\la \tau\ra}$ are the fixed
point subspaces of $\tau$ in $V$ and $V^T$ respectively. We shall
show that the weight one subspace $V(\tau)_1$ of $V(\tau)$ is
nontrivial and the Lie algebra structure on $V(\tau)_1$ is abelian.
Hence, by a result of Dong and Mason \cite[Theorem 3]{dm2},
$V(\tau)$ is isomorphic to the Leech lattice VOA $V_\Lambda$.
Therefore, by reversing the orbifold construction on $V(\tau)$, we
obtain that $V^{\la\tau\ra}\cong V_\Lambda^+$ and $V\cong
V_\Lambda^+\oplus V_\Lambda^{T,+}\cong V^\natural$ as
$V_\Lambda^+$-module. It is well-known that there exists a unique
simple vertex operator algebra structure on $V_\Lambda^+\oplus
V_\Lambda^{T,+}$ (cf.\ \cite{dm1,huang}) and thus we conclude that
$V$ and $V^\natural$ are isomorphic simple vertex operator algebras.

Another uniqueness result of the moonshine VOA $V^\natural$ is also
obtained in \cite{dgl} under a different set of assumptions. Again,
the existence of $48$ commuting Virasoro algebras of central charge
$1/2$ is crucial to their argument.

\medskip

\paragraph{\bf Acknowledgment}
Part of the work was done when the second author was visiting the
National Center for Theoretical Sciences, Taiwan on August 2006. He
thanks the staffs of the center for their help.

\section{Framed vertex operator algebra}
First we shall recall the definition of framed vertex operator
algebras and review several important results \cite{dgh,ly,m2,m3,m4}.

\begin{df}\label{df:3.2}
  (\cite{dgh,m3})
  A simple vertex operator algebra $(V,\om)$ is called \emph{framed} if there exists a
  mutually orthogonal set  $\{e^1, \dots,e^n\}$ of Virasoro elements of central charge $1/2$
  such that $\om=e^1+\cds +e^n$ and each $e^i$ generates a simple Virasoro vertex operator algebra
  isomorphic to $L(\shf,0)$.
  The full sub VOA  $F$ generated by $e^1,\dots,e^n$ is called an \emph{Virasoro frame} or simply
  a \emph{frame} of $V$.
\end{df}

Let $(V,\om)$ be a framed VOA with a frame $F$.
We denote by $\vir(e^i)$ the Virasoro subVOA generated by $e^i$.
Then
$$
  F \cong \vir(e^1)\tensor \cds \tensor \vir(e^i)\cong  L(\shf,0)^{\tensor
  n}.
$$
Since the Virasoro VOA $L(\shf,0)$ is rational, $V$ is a direct sum
of irreducible $F$-submodules $\tensor_{i=1}^n L(\shf,h_i)$ with
$h_i\in \{ 0,1/2,1/16\}$, namely
\[
  V= \bigoplus_{h_i\in \{ 0,1/2, 1/16\}} m_{h_1,\dots,h_n}
       L(\shf,h_1)\tensor \cds \tensor L(\shf,h_n),
\]
where $m_{h_1,\dots,h_n}\in \N$ denotes the multiplicity.

For each irreducible $F$-module
$\tensor_{i=1}^n L(\shf,h_i)$, we define its binary {\it
$1/16$-word} (or {\it $\tau$-word}) $(\alpha_1,\cds,\alpha_n)\in
\Z_2^n$ by $\alpha_i=1$ if and only if $h_i=1/16$. For $\alpha \in
\Z_2^n$, denote by $V^\alpha$ the sum of all irreducible
$F$-submodules of $V$ whose $1/16$-words are equal to $\alpha$.
Define $D=\{ \alpha \in \Z_2^n \mid V^\alpha\ne 0\}.$ Then $D$ is a
linear code and  we have the {\it 1/16-word decomposition}
$$
  V=\bigoplus_{\alpha \in D} V^\alpha.
$$
It is shown in Dong et al.\ \cite{dmz} that
\[
  V^0= \bigoplus_{h_i\in \{ 0,1/2\}} m_{h_1,\dots,h_n}
       L(\shf,h_1)\tensor \cds \tensor L(\shf,h_n)
\]
is a subalgebra of $V$ and the multiplicity $m_{h_1,\cds,h_n}\leq 1$
for $h_i\in \{0,1/2\}.$
Thus we obtain another linear code $C:= \{
(2h_1,\cds ,2h_n)\in \Z_2^n \mid h_i\in \{ 0,\shf\},\
m_{h_1,\cds,h_n}\ne 0\}$  and $V^0$ can be decomposed as
\begin{equation}\label{eq:3.2}
  V^0=\bigoplus_{\alpha =(\alpha_1,\cds,\alpha_n)\in C}
  L(\shf,\alpha_1/2) \tensor \cds \tensor L(\shf,\alpha_n/2).
\end{equation}
The VOA $V^0$ is often called the code vertex operator algebra
associated with the code $C$ and denoted by $M_C$. Note that the VOA
structure of $V^0$ is uniquely determined by the code $C$
\cite{dgh,m2}.

Summarizing, there exists a pair  of even linear codes $(C,D)$ such
that $V$ is an $D$-graded extension of a code VOA $M_C$ associated
to $C$. We call the  pair $(C,D)$ the {\it structure codes} of a
framed VOA $V$ associated with the frame $F$. Since the powers of
$z$ in an $L(\shf,0)$-intertwining operator of type
$L(\shf,\shf)\times L(\shf,\shf)\to L(\shf,\sfr{1}{16})$ are
half-integral, the structure codes $(C, D)$ satisfy $C\subset
D^\perp$. Moreover, $C=D^\perp$ if and only if $V$ is holomorphic
(cf.\ \cite{dgh,ly,m4}).


Let $V$ be a framed VOA  with the structure codes $(C,D)$, i.e.,
$V=\oplus_{\al\in D} V^\al$, and $V^0=M_C$. For a binary codeword
$\beta \in \Z_2^n$, we define a linear map $\tau_\beta: V\to V$ by
\begin{equation}\label{eq:3.3}
 \tau_\beta=(-1)^{\la \al,\be\ra} v, \quad \text{ for } v\in V^\al.
\end{equation}
By the fusion rules, it is easy to show that $\tau_\beta$ is an
automorphism of $V$ (cf.\ \cite{m1}). This automorphism is often
called a $\tau$-involution (or Miyamoto involution).
Let $P=\{\tau_\beta \mid \beta\in \Z_2^n\}$ be the subgroup generated by the
$\tau$-involutions. Then $P\cong \Z_2^n/ D^\perp$ and $V^0=V^P$ is
the fixed point subalgebra. Thus, all $V^\al, \al\in D,$ are
irreducible modules over $V^0=M_C$ (cf.\ \cite{DLM2}).

\medskip

Next we shall recall a very important result from \cite{ly}.
\begin{thm}[Theorem 5.6 of \cite{ly}]\label{thm:5.6}
  Let $V=\oplus_{\alpha\in D} V^\alpha$ be a framed VOA with structure codes $(C,D)$.
  Then $V^\al, \al\in D,$ are all simple current modules over
  the code VOA $V^0=M_C$.
\end{thm}

The following two results follow immediately from the general
arguments on simple current extensions \cite{dlm,lam2,y2}.

\begin{thm}[ Corollary 5.7 of \cite{ly}]\label{cor:5.7}
  Let $V=\oplus_{\alpha\in D} V^\alpha$ be a framed VOA with structure codes $(C,D)$.
  Let $W$ be an irreducible $V^0$-module.
  Then there exists a unique $\eta\in \Z_2^n$ up to $D^\perp$ such that
  $W$ can be uniquely extended to an irreducible $\tau_\eta$-twisted $V$-module
  which is given by $V \fusion_{V^0} W$ as a $V^0$-module.
\end{thm}

\begin{thm}[cf.\ \cite{lam3,ly}]\label{2.4}
Let $V=\oplus_{\alpha\in D} V^\alpha$ be a holomorphic framed VOA
with structure codes $(C,D)$. For any $\delta\in \Z_2^n$, denote
\[
D^0=\{\al\in D\,|\, \la \al, \delta\ra=0\}\quad \text{ and }\quad
D^1=\{\al\in D\,|\, \la \al, \delta\ra\neq 0\}.
\]
Define
\[
V(\tau_\delta)=
\begin{cases}\ds
(\bigoplus_{\al\in D^0} V^\al) \oplus (\bigoplus_{\al\in D^1}
M_{\delta+C}\fusion_{M_C} V^\al) & \text{ if } \wt\,\delta \text{ is
odd},\\
\ds (\bigoplus_{\al\in D^0} V^\al) \oplus (\bigoplus_{\al\in D^0}
M_{\delta+C}\fusion_{M_C} V^\al) & \text{ if } \wt\,\delta \text{ is even}.
\end{cases}
\]
Then $V(\tau_\delta)$ is also a holomorphic framed VOA. Moreover,
the structure codes of $V(\tau_\delta)$ are given by $(C,D)$ if
$\wt\,\delta$ is odd and $(C\cup(\delta+C), D^0)$ if $\wt\,\delta$
is even.
\end{thm}

\begin{rem}
The above construction of the holomorphic VOA $V(\tau_\delta)$ is
referred to as a $\tau_\delta$-orbifold construction of $V$.
\end{rem}

\section{Strongly rational, holomorphic vertex operator algebra}
In this section, we shall review some basic facts about strongly
rational vertex operator algebra from \cite{dm1,dm2}.

\begin{df}
A VOA $V$  is of \emph{CFT-type} if the natural $\Z$-grading on $V$ takes
the form $V = V_0 \oplus V_1 \oplus\cdots$ with $V_0 = \C\vac$.
\end{df}

\begin{rem}
Let $V=\oplus_{n=0}^\infty V_n$ be a VOA of CFT-type. Then the
weight one subspace $V_1$ carries a structure of a Lie algebra with the
bracket
\[
[u,v]=u_0v,\quad u,v\in V_1
\]
and an invariant bilinear form defined by
\[
\la u,v\ra \vac = u_1v\quad \text{ for } u,v\in V_1.
\]
\end{rem}

\begin{df}[\cite{dm1}]
A vertex operator algebra $V$ is said to be {\it strongly
rational} if it satisfies the following conditions:\\
1. $V$ is of  CFT type and $L(1)V_1=0$.\\
2. V is $C_2$-cofinite, i.e., $\dim V/C_2(V) <\infty$,
   where $C_2(V)=\mathrm{span}\{ u_{-2}v \mid u,v\in V\}$. \\
3. V is rational.
\end{df}

\begin{rem}
All framed vertex operator algebras are strongly rational.
\end{rem}

Recently, Dong and Mason\,\cite{dm1,dm2,dm3} have a study of
strongly rational, holomorphic vertex operator algebras of rank
$\leq 24$ by using the Lie algebra structure on $V_1$. Although their method
is not very effective when $V_1=0$, they obtained the following
characterization of the Leech lattice VOA $V_\Lambda$.

\begin{thm}[Dong-Mason\,\cite{dm2}]\label{leech}
Let $V$ be a strongly rational, holomorphic vertex operator algebra
of rank $24$ such that the Lie algebra on $V_1$  is abelian. Then
$\dim V_1=24$ and $V$ is isomorphic to the Leech lattice VOA
$V_\Lambda$.
\end{thm}

\section{Uniqueness of the moonshine vertex operator algebra}

Let $V$ be a holomorphic framed VOA of rank $24$ such that $V_1=0$.
We shall prove that $V$ is isomorphic to the moonshine VOA
$V^\natural$ by using the theorem of Dong-Mason (cf.\ Theorem \ref{leech}).
First let $C$ and $D$ be the structure codes of $V$.
That means $$V=\oplus_{\be\in D} V^\be, \quad V^0=M_C\quad \text{
and }\quad  C=D^\perp.$$
Since $V_1=0$, the code $C$ contains no
element of weight $2$. Now let $\delta= (110\dots 0)$. Then
$\delta\notin C=D^\perp$ and hence the automorphism
$\tau_{\delta}$ defines a (nontrivial) involution on $V$ and the
fixed point subVOA is given by
\[
V^{\la \tau_{\delta}\ra}= \bigoplus_{\be\in D^0} V^\be, \quad \text{
where } D^0=\{ \be \in D \mid  \la \be, \delta\ra =0\} .
\]

Now define the $\tau_\delta$-orbifold construction of $V$ by
\[
V(\tau_\delta) = \bigoplus_{\be\in D^0} \left( V^\be \oplus
M_{\delta+C}\fusion_{M_C} V^\be \right)
\]

By Theorem \ref{2.4}, $V(\tau_{\delta})$ is a holomorphic framed VOA
of rank $24$ and the structure codes of $V(\tau_\delta)$ is given by
$(\tilde{C}, D^0)$, where $\tilde{C}$ is the binary code generated by
$C$ and $\delta$.

\begin{rem}\label{g}
One can define an involution $g$ on $V(\tau_\delta)$ as
follows:
\[
g=
\begin{cases}
\ \ 1& \text{ on } \bigoplus_{\be\in D^0} V^\be,\\
-1& \text{ on } \bigoplus_{\be\in D^0} M_{\delta+C}\fusion_{M_C}  V^\be.\\
\end{cases}
\]
Moreover, the fixed point subalgebra
$V(\tau_\delta)^{\la g\ra}=V^{\la \tau_\delta\ra}=\bigoplus_{\be\in D^0} V^\be$.
\end{rem}

\begin{prop}
$V(\tau_{\delta})_1\neq 0$ and the Lie algebra on
$V(\tau_{\delta})_1$ is abelian. As a consequence,
$V(\tau_{\delta})$ is isomorphic to the Leech lattice VOA
$V_\Lambda$.
\end{prop}

\begin{proof}
Since $M_{\delta+C}\subset V(\tau_{\delta})$ and $\wt\,\delta=2$,
we have $V(\tau_{\delta})_1\neq 0$. Now let $g$ be the involution
on $V(\tau_{\delta})$ defined in Remark \ref{g}. Denote by
$V(\tau_{\delta})^+$ and $V(\tau_{\delta})^-$ the fixed point
subspace and the $-1$-eigenspace of $g$ on $V(\tau_{\delta})$,
respectively. Since $ V(\tau_{\delta})^+= V^{\la \tau_{\delta}\ra}\subset
V$, the weight one subspace $V(\tau_{\delta})^+_1$ of $V(\tau_\delta)$ is trivial
and hence $V(\tau_{\delta})_1= V(\tau_{\delta})^-_1$.
Therefore, for any $u,v\in V(\tau_{\delta})_1$,
we have $u,v\in V(\tau_{\delta})^-$ and thus
$$[u,v]=u_0v\in [V(\tau_\delta)^-_1,V(\tau_\delta)^-_1]
\subset V(\tau_{\delta})^+_1=0.$$
Therefore, the Lie algebra
$V(\tau_\delta)_1$ is abelian, and hence $V(\tau_\delta)$ is
isomorphic to the Leech lattice VOA $V_\Lambda$ by Theorem
\ref{leech}.
\end{proof}

\begin{thm}
Let $V$ be a holomorphic framed VOA of rank $24$ such that $V_1=0$.
Then  $V$ is isomorphic to the moonshine VOA $V^\natural$.
\end{thm}

\begin{proof}
By the above proposition, we know that $V(\tau_{\delta})$ is
isomorphic to the Leech lattice VOA $V_\Lambda$.
Moreover, $g$ acts on $V(\tau_{\delta})_1$ as $-1$ so that $g$ is conjugate
to the lift $\theta$ of $-1$ map on $\Lambda$ (cf.\ Eq.\ 1.1) by Theorem D.6 of
\cite{dgh} (see also Lemma 7.15 of \cite{ly}).
Therefore, we have $V(\tau_\delta)^{\la g\ra}=V^{\la \tau_{\delta}\ra}\cong V_\Lambda^+$.
It is well-known (cf.\ \cite{d1,huang}) that $V_\Lambda^+$ has exactly $4$
inequivalent irreducible modules, namely $V_\Lambda^+$, $V_\Lambda^-$,
$V_\Lambda^{T,+}$ and $V_\Lambda^{T,-}$, and all these modules are simple currents.
Their top weights are $0, 1, 2$ and $3/2$, respectively.
Since $V$ has integral weights and $V_1=0$, it is clear that
\[
V\cong V_\Lambda^+ \oplus V_\Lambda^{T,+}\cong V^\natural
\]
as $V_\Lambda^+$-modules.
Then by the uniqueness of simple current extensions shown in \cite{dm1},
we can establish the desired isomorphism between $V$ and $V^\natural$.
\end{proof}

\bibliographystyle{amsplain}

\end{document}